\documentclass{amsart}

\newtheorem{theorem}{Theorem}[section]
\newtheorem{lemma}[theorem]{Lemma}
\newtheorem{corollary}[theorem]{Corollary}

\theoremstyle{definition}

\theoremstyle{remark}
\newtheorem{remark}[theorem]{Remark}

\numberwithin{equation}{section}

\usepackage{hyperref}

\begin{document}

\title[Some nonlinear inequalities]
{Some nonlinear inequalities and applications}

\author{N. S. Hoang$\dag$}
\address{Mathematics Department, Kansas State University,
Manhattan, KS 66506-2602, USA}
\email{nguyenhs@math.ksu.edu}

\author{A. G. Ramm$\ddag$}%\footnotemark[1]}
\address{Mathematics Department, Kansas State University,
Manhattan, KS 66506-2602, USA}
\email{ramm@math.ksu.edu}

%    General info
%\renewcommand{\thefootnote}{\fnsymbol{footnote}}
%\footnotetext[1]{Corresponding author}
\date{}

\begin{abstract}
Sufficient conditions are given for the relation $\lim_{t\to\infty}y(t)=0$ 
to hold, 
where $y(t)$ is a continuous nonnegative function on $[0,\infty)$ 
 satisfying some nonlinear inequalities.
The results are used for a study of large time behavior of 
the solutions to nonlinear 
evolution equations. 
Example of application is given for a solution to some evolution 
equation with a nonlinear partial differential operator.

\end{abstract}


\maketitle

\section{Introduction}

The stability study of many evolution equations is 
 a study of  large time  behavior of the solutions to these equations. 
In this paper we reduce such a study to a study of the behavior of a 
solution  $y(t)$ to some nonlinear inequalities. Assume that a nonnegative 
continuous function $y(t)$ satisfies the following conditions 
\begin{equation}
\label{eq1}
\int_0^\infty \omega(y(t))\frac{1}{(1+t)^\alpha}dt<\infty,
\qquad 0\le \alpha\le 1,
\end{equation}
and 
\begin{equation}
\label{eq2}
y(t) - y(s) \le \int_s^t f(x,y(x))dx,\qquad 0\le s\le t,
\end{equation}
where $f(x,y)$ is a nonnegative continuous function on $[0,\infty)\times 
[0,\infty)$,   
$0\le\omega(t)$ is a non-decreasing continuous function, and   
$\omega(t)=0$ implies $t=0$. 

The question arises: 

{\it Under what condition on $f(y,t)$ does it follow that \begin{equation}
\label{eq3} \lim_{t\to\infty} y(t) = 0? \end{equation} } 

There is a very large literature on inequalities (see, e.g.,  \cite{Be}, 
\cite{Bu}
and references therein). The Barbalat's lemma is 
an integral inequality used in applied nonlinear control (\cite{S}).
The inequalities, derived in this
paper,  are new and are useful in many applications. 
In
\cite[p.227]{zheng} inequality \eqref{eq1} is studied for $\omega(t) = t$
and $\alpha=0$.  In this case condition \eqref{eq1} becomes $y\in
L^1[0,\infty)$.  In \cite{zheng} it is proved that
\eqref{eq3} holds if  $y\in
L^1[0,\infty)$ and the following two conditions hold: 
\begin{equation} \label{eq4} y(t) - y(s) \le
\int_s^t f(y(x))dx + \int_s^t h(x)dx,\qquad \int_0^\infty h(x)dx<\infty.  
\end{equation} Here $f$ and $h$ are nonnegative functions, and $f$ is
continuous  and non-decreasing.  Proofs of this result can be found in
\cite{zheng} and in \cite{R499}.  Applications of this result 
to the stability study of evolution equations can be found in 
\cite{zheng} and references therein.
This result is not
applicable if $y(t) = O(\frac{1}{t^\beta})$ as $t\to \infty$,
where $\beta\in(0,1)$, because then $y(t)$ is not in $L^1[0,\infty)$.  
Also, this result is not applicable if \eqref{eq2} holds instead of 
\eqref{eq4} and
$f$ depends on $x$.

The second nonlinear inequality we study is the following one:
\begin{equation}
\label{eq**}
\dot{g}(t) \le -a(t)f(g(t)) + b(t),\qquad t\ge 0,
\end{equation}
where $a,b$ and $g$ are nonnegative functions on $[0,\infty)$,
$g\in C^1([0,\infty)),$ $a\in C([0,\infty))$ and $b\in 
L^1_{loc}([0,\infty))$. 
A sufficient condition for the relation 
$\lim_{t\to\infty}g(t)=0$ to hold is proposed and justified in 
\cite{R499}.
In our paper inequality 
\eqref{eq**} is studied by a different method and some new
sufficient conditions for \eqref{eq3} to hold are proposed and justified.

The paper is organized as follows. In Theorems 2.1, 2.4 and 
\ref{theorem2} and their corollaries, sufficient conditions for
\eqref{eq4} to hold are formulated and justified. In Theorems 2.11, 2.13
and 2.14 sufficient conditions
for the relation $\lim_{t\to\infty}g(t)=0$ to hold are proposed and
justified under the assumption that $f(t)$ is a continuous and
non-decreasing function on $[0,\infty)$.  In Section 3 applications of the 
new
results to the stability study of evolution equations are given.

\section{Main results}

Throughout the paper we assume that 

{\it $\omega(t)\geq 0$ is 
a non-decreasing continuous function  and if $\omega(t)=0$ then $t=0$.} 

This assumption is standing and is
not repeated. 

\begin{theorem}
\label{maintheorem}
Let $y(t)\ge 0$ be a continuous function on $[0,\infty)$,
\begin{equation}
\label{eq1.1}
\int_0^\infty \omega(y(t))dt <\infty,
\end{equation}
and
\begin{equation}
\label{eq5}
y(t) - y(s) \le \int_s^t f(\xi,y(\xi))d\xi,\qquad 0\le s\le t,
\end{equation}
where $f(t,y)$ is a nonnegative continuous function on 
$[0,\infty)\times [0,\infty)$. 
Define
\begin{equation}
\label{eq6}
F(t,v) := \int_0^t \max_{0\le \zeta\le v} f(\xi,\zeta)d\xi,\qquad v,t \ge 
0.
\end{equation}
If there exists a constant $a>0$ such that the function $F(t,a)$ is 
uniformly continuous with respect to $t$ on $[0,\infty)$,  
then
\begin{equation}
\label{eq7}
\lim_{t\to\infty}y(t) = 0.
\end{equation}
\end{theorem}

\begin{proof}
If \eqref{eq7} does not hold, then there exists an $\epsilon>0$ 
and a sequence $(t_n)_{n=1}^\infty$ such that
\begin{equation}
\label{eq8}
0<t_n \nearrow \infty,\qquad y(t_n) \ge \epsilon,\qquad \forall n\ge 1. 
\end{equation}
Without loss of generality we assume that $\epsilon<a$. 

Since $F(t,a)$ is uniformly continuous with respect to $t$, there exists 
$\delta>0$ 
such that 
\begin{equation}
\label{eq9}
\int_{t}^{t+\delta} \max_{0\le y\le a} f(\xi,y)d\xi 
= F(t+\delta,a) - F(t,a) < \frac{\epsilon}{2},\qquad \forall t\ge 0.
\end{equation}

Let us prove that 
\begin{equation}
\label{eq10}
y(t) \ge \frac{\epsilon}{2},\qquad \forall t\in [t_n-\delta, t_n],
\qquad \forall n\geq 1.
\end{equation}
Assume that \eqref{eq10} does not hold. Then there exists 
$\tilde{n}>0$ and $\xi \in  [t_{\tilde{n}}-\delta, t_{\tilde{n}})$ 
such that 
\begin{equation}
\label{eq11}
y(\xi)< \frac{\epsilon}{2}.
\end{equation}
Let 
\begin{equation}
\label{eq12}
\nu = \min \{x: \xi< x\le t_{\tilde{n}},\, y(x) \ge \epsilon\}.
\end{equation}
From the continuity of $y$, \eqref{eq8}, and \eqref{eq11}--\eqref{eq12} 
one obtains 
\begin{equation}
\label{eq13}
t_{\tilde{n}} - \delta \le \xi < \nu \le t_{\tilde{n}},\qquad 
y(\nu)=\epsilon,
\end{equation}
and
\begin{equation}
\label{eq14}
0 \le y(x) \le y(\nu) = \epsilon,\qquad \xi\le x\le \nu.
\end{equation}
It follows from \eqref{eq5}, \eqref{eq11}, and \eqref{eq13}--\eqref{eq14} that
\begin{equation}
\label{eq15}
\begin{split}
\frac{\epsilon}{2} < y(\nu) - y(\xi) 
&\le \int_\xi^\nu f(x,y(x))dx \le \int_\xi^\nu \sup_{0 \le \zeta\le 
\epsilon} f(x,\zeta)dx \\
&\le \int_{t_n-\delta}^{t_n} \sup_{0 \le \zeta\le a} f(x,\zeta)dx 
< \frac{\epsilon}{2}.
\end{split}
\end{equation}
This contradiction proves that \eqref{eq10} holds. 

From \eqref{eq10} one gets
\begin{equation}
\label{eq16}
\int_{t_n - \delta}^{t_n} \omega(y(x))dx \ge \delta 
\omega(\frac{\epsilon}{2})>0,\qquad \forall n \ge 1.
\end{equation}
This contradicts the Cauchy criterion for the convergence
of the integral (2.1).% $\omega(y(t))\in L^1[0,\infty)$. 
Thus,  \eqref{eq7} holds. 

Theorem~\ref{maintheorem} is proved.
\end{proof}

\begin{remark}{\rm
If $F(t,a)$ is uniformly continuous with respect to $t$ on 
$[0,\infty)$, then $F(t,v)$ is uniformly continuous with respect to 
$t$ on $[0,\infty)$ for all $v\in [0,a]$. However, $F(t,v)$ may be not 
uniformly continuous with respect to $t$ 
on $[0,\infty)$ for some $v>a$. Here is an example: 

Let
\begin{equation}
\label{eqju1}
f(x,y) := \left \{\begin{matrix} 1&\quad \text{if}\quad  & 0\le y\le 1\\
1 + (y-1)t & \quad \text{if}\quad  & y\ge 1.
\end{matrix}
\right. 
\end{equation}
By a simple calculation one gets
\begin{equation}
\label{eqju2}
F(t,u) = \left \{\begin{matrix} t&\quad \text{if}\quad & 0\le u\le 1\\
t + (u-1)\frac{t^2}{2} &\quad \text{if}\quad & u\ge 1.
\end{matrix}
\right. 
\end{equation} 
It follows from \eqref{eqju2} that $F(t,u)$ is uniformly 
continuous with respect to $t$ on $[0,\infty)$ if and only if $u\in [0,1]$. 
}
\end{remark}

From Theorem~\ref{maintheorem} one derives the following corollary. 
\begin{corollary}
\label{corollary2.2}
Assume that $y(t)\geq 0$ be a  continuous function 
satisfying inequality \eqref{eq1.1}, 
\begin{equation}
\label{eq18}
y(t) - y(s) \le \int_s^t [g(\xi)\varphi(y(\xi)) + h(\xi)]d\xi,
\qquad 0\le s\le t,
\end{equation}
where $g$ and $h$ are nonnegative locally integrable functions on 
$[0,\infty)$,  
$\varphi\ge 0$ is a continuous function on $[0,\infty)$, and 
the functions $\int_0^tg(x)dx$ and 
$\int_0^t h(x)dx$ are uniformly continuous with respect to $t$ on $[0,\infty)$. 
Then \eqref{eq7} holds.
%\begin{equation}
%\label{eq19}
%\lim_{t\to\infty}y(t) = 0.
%\end{equation}
\end{corollary}

\begin{proof}
Let 
$$
f(x,y) := g(x)\varphi(y) + h(x),\qquad x\ge 0,\qquad y\ge 0.
$$
It follows from the uniform continuity of $\int_0^tg(x)dx$ and 
$\int_0^t h(x)dx$ that the function $F$, defined in \eqref{eq6}, 
is uniformly continuous. Thus, \eqref{eq7} follows 
from Theorem~\ref{maintheorem}. 
\end{proof}

\begin{theorem}
\label{theorem2.01}
Assume that $y(t)\ge 0$ is a continuous function on $[0,\infty)$,   
\begin{equation}
\label{teq20}
\int_0^\infty \omega(y(t))\varphi(t)dt <\infty,
\end{equation}
where $\varphi(t)\geq 0$ is a continuous function on $[0,\infty)$, and 
there exists a constant $C>0$ such that
\begin{equation}
\label{teq20.0}
\lim_{t\to\infty} \Big(t - \frac{C}{\varphi(t)}\Big) = \infty,\qquad 
M:=\limsup_{t\to\infty}\frac{\max_{\xi\in[t-\frac{C}{\varphi(t)},t]}
\varphi(\xi)}
{\min_{\xi\in[t-\frac{C}{\varphi(t)},t]}\varphi(\xi)} <\infty,
\end{equation}
where  $f(x,y)\geq 0$ is a continuous on $[0,\infty)\times 
[0,\infty)$ function, which satisfies condition (2.2).
%\begin{equation}
%\label{teq21}
%y(t) - y(s) \le \int_s^t f(\xi,y(\xi))d\xi,\qquad 0\le s\le t.
%\end{equation}
If there exist  constants $a>0$ and $\theta >0$ such 
that the following condition holds:
\begin{equation}
\label{teq22}
\int_s^t \sup_{0\le \zeta\le a}f(x,\zeta)dx\le  
(t-s)\theta a\max_{\xi\in[s,t]}\varphi(\xi),
\qquad \theta=const >0, 
\quad t>s\gg 1,
\end{equation}
then \eqref{eq7} holds.
%\begin{equation}
%\label{teq24}
%\lim_{t\to\infty}y(t) = 0.
%\end{equation}
\end{theorem}

\begin{remark}{\rm
In \eqref{teq22} and below the notation $s\gg 1$ means "for all 
sufficiently large $s>0$". 
}
\end{remark}

\begin{proof}
{\it Let us consider first Case 1, namely,  $0<\theta <1$. Later we 
reduce Case 2, namely, $\theta\geq 1$, to Case 1.} 

Assume that \eqref{eq7} does not hold. Then there exists an 
$\epsilon>0$, 
a sequence $(t_n)_{n=1}^\infty$ such that 
\begin{equation}
\label{tzzeq8}
0 < t_n \nearrow \infty,\qquad y(t_n) \ge \epsilon,\qquad \forall n\ge 1,
\end{equation}
and without loss of generality one assumes that 
\begin{equation}
\label{teqq02}
\epsilon \le 2a MC. 
\end{equation}

Let us prove that 
\begin{equation}
\label{tzzeq10}
y(t) \ge (1-\theta)\epsilon,\qquad \forall t\in 
[\tilde{t}_n, t_n],\qquad \forall n\gg 1,
\end{equation}
where
\begin{equation}
\label{eqfff5}
\tilde{t}_n:= t_n - \frac{\epsilon}{2aM\varphi(t_n)}<t_n.
\end{equation}
Assume that \eqref{tzzeq10} does not hold. Then there exists 
a  sufficiently large $\tilde{n}>0$ and a 
$\xi \in  [\tilde{t}_{\tilde{n}}, t_{\tilde{n}})$ 
such that 
\begin{equation}
\label{tzzeq11}
y(\xi)< (1-\theta)\epsilon.
\end{equation}
Let 
\begin{equation}
\label{tzzeq12}
\nu = \min \{x: \xi\le x\le t_{\tilde{n}},\, y(x) \ge \epsilon\}.
\end{equation}
Then 
\begin{equation}
\label{tzzeq13}
\tilde{t}_{\tilde{n}} \le \xi < \nu \le t_{\tilde{n}},
\end{equation}
and
\begin{equation}
\label{tzzeq14}
0 \le y(x) \le y(\nu) = \epsilon,\qquad \xi\le x\le \nu.
\end{equation}
It follows from (2.2), \eqref{tzzeq8}, \eqref{tzzeq11}, 
and \eqref{tzzeq13}--\eqref{tzzeq14} that
\begin{equation}
\label{tzzeq15}
\begin{split}
\theta\epsilon < y(\nu) - y(\xi) 
&\le \int_\xi^\nu f(x,y(x))dx \le \int_\xi^\nu \sup_{0 \le \zeta\le 
\epsilon} f(x,\zeta)dx \\
&\le \int_{\tilde{t}_{\tilde{n}}}^{t_{\tilde{n}}} \sup_{0 \le \zeta\le a} f(x,\zeta)dx \le 
(t_{\tilde{n}} - \tilde{t}_{\tilde{n}})\theta a
\max_{\xi\in[\tilde{t}_{\tilde{n}},t_{\tilde{n}}]}\varphi(\xi)\\
&=\theta a \frac{\epsilon \max_{\xi\in[\tilde{t}_{\tilde{n}},t_{\tilde{n}}]}
\varphi(\xi)}{2Ma\varphi(t_{\tilde{n}})} \le \theta \epsilon.
\end{split}
\end{equation}
This contradiction proves \eqref{tzzeq10}. 
In the derivation of \eqref{tzzeq15} we have used the following 
inequality: 
\begin{equation}
\label{eqiu3}
\frac{\max_{\xi\in[\tilde{t}_{\tilde{n}},t_{\tilde{n}}]}\varphi(\xi)}{\varphi(t_{\tilde{n}})} 
\le \frac{\max_{\xi\in[t_{\tilde{n}} - C\varphi^{-1}(t_{\tilde{n}}), t_{\tilde{n}}]}\varphi(\xi)}
{\min_{\xi\in[t_{\tilde{n}} - C\varphi^{-1}(t_{\tilde{n}}), t_{\tilde{n}}]}\varphi(\xi)} 
< 2M,\qquad \tilde{n} \gg 1,
\end{equation}
which follows from \eqref{teq20.0} for sufficiently large $t_{\tilde{n}}$,
and the factor $2$ in \eqref{eqiu3} can be replaced by any fixed factor 
$1+q$, where
$q>0$ can be arbitrarily small if $t_{\tilde{n}}$ is sufficiently large. 

Since $\omega(t)$ is non-decreasing, it follows from \eqref{tzzeq10} that 
\begin{equation}
\label{tzzeq16}
\begin{split}
\int_{\tilde{t}_n}^{t_n} \omega(y(x))\varphi(x)dx 
& \ge  (t_n - \tilde{t}_n) 
\omega\big{(}(1-\theta)\epsilon\big{)}\min_{\tilde{t}_n\le \xi\le 
t_n}\varphi(\xi)\\
%%& \ge \omega\big{(}(1-\theta)\epsilon\big{)}\frac{\epsilon\varphi(t_n)}{aM\varphi(t_n)}\\
& \ge \omega\big{(}(1-\theta)\epsilon\big{)} \frac{\epsilon}{2aM}
\frac{\min_{\tilde{t}_n\le\xi\le t_n}\varphi(\xi)}{\max_{\tilde{t}_n\le\xi\le t_n}\varphi(\xi)}\\
&\ge  \omega\big{(}(1-\theta)\epsilon\big{)} \frac{\epsilon}{2aM(M+q)}>0,
\end{split}
\end{equation}
where $q>0$ is arbitrarily small for all  sufficiently large $n$. 
From \eqref{teqq02}, \eqref{teq20.0}, and \eqref{tzzeq8},  one gets
\begin{equation}
\label{teqq01}
\lim_{n\to\infty}\bigg{(} t_n - \frac{\epsilon}{2a M\varphi(t_n)} \bigg{)} 
\ge \lim_{n\to\infty}\bigg{(} t_n - \frac{C}{\varphi(t_n)} \bigg{)}
=\infty.
\end{equation}

Inequalities \eqref{tzzeq16} and \eqref{teqq01} contradict the 
Cauchy criterion for the convergence of integral \eqref{teq20}. 
Thus, \eqref{eq7} holds. 

{\it Consider Case 2, namely  $\theta\ge 1$}. 
In this case one replaces $\theta$ by $\theta_1=\frac 1 2,$ $C$ by 
$C_1=2\theta C$, $M$ by $M_1=M$, defined
in \eqref{teq20.0} with the $C_1$ in place of $C$, and, therefore, one 
reduces 
the problem to Case 1 with $\theta=\frac 1 2<1$.

Let us give a more detailed argument.
Let $\varphi_1(t) := 2\theta\varphi(t)$ 
and $C_1: = 2\theta C$. Then 
\begin{equation}
\label{ueq0}
\frac{C_1}{\varphi_1(t)} = \frac{C}{\varphi(t)},\qquad \forall t\ge 0.
\end{equation}
This, \eqref{teq20}, \eqref{teq20.0} and \eqref{teq22} imply
\begin{gather}
\label{ueq1}
\int_0^\infty \omega(y(t))\varphi_1(t)dt <\infty,\\
\label{ueq2}
\lim_{t\to\infty} \Big(t - \frac{C_1}{\varphi_1(t)}\Big) = \infty,\qquad 
\limsup_{t\to\infty}\frac{\max_{\xi\in[t-\frac{C_1}{\varphi_1(t)},t]}
\varphi_1(\xi)}
{\min_{\xi\in[t-\frac{C_1}{\varphi_1(t)},t]}\varphi_1(\xi)} = M<\infty\\
\label{ueq3}
\int_s^t \sup_{0\le \zeta\le a}f(x,\zeta)dx\le  
(t-s)\frac{a}{2}\max_{\xi\in[s,t]}\varphi_1(\xi),
\qquad t>s\gg 1.
\end{gather}
%From the proof above one gets \eqref{teq24}. 
Theorem~\ref{theorem2.01} is proved.
\end{proof}

\begin{remark}
\begin{enumerate}
\item[(i)]{Conditions \eqref{teq20.0} hold if 
\begin{equation}
\label{equi8}
\liminf_{t\to\infty} t\varphi(t) > 0,\qquad \limsup_{t\to\infty} 
\frac{\max_{\xi\in[(1-\epsilon)t, t]}\varphi(\xi)}{\min_{\xi\in[(1-\epsilon)t, t]}\varphi(\xi)} < \infty,
\end{equation}
for a sufficiently small $\epsilon>0$.} 

\item[(ii)]{
If $y(t)$ is differentiable, then (2.2) is equivalent to
\begin{equation}
\label{eqv1}
y'(t) \le f(t,y(t)),\qquad t\ge0.
\end{equation}
}
\item[(iii)]{
Theorem~\ref{theorem2.01} holds if in place of  \eqref{teq22} one assumes 
that 
\begin{equation}
\label{eqv2}
\sup_{0\le\zeta\le a}f(t,\zeta) \le \tilde{C} \varphi(t),\qquad t\gg 1,\qquad \tilde{C}=const>0.
\end{equation}
Indeed, if \eqref{eqv2} hold then
$$
\int_s^t \sup_{0\le\zeta\le a}f(\xi,\zeta)d\xi \le
\int_s^t \tilde{C} \varphi(\xi)d\xi 
\le \tilde{C}(t-s)\max_{s\le \xi\le t} \varphi(\xi). 
$$
}
\item[(iv)]{
If $\varphi(t)$ is non-increasing, then the second relation in 
\eqref{teq20.0} becomes 
\begin{equation}
\label{eqv4}
M:=\limsup_{t\to\infty}\frac{\varphi(t-\frac{C}{\varphi(t)})}
{\varphi(t)} < \infty.
\end{equation}
}
\end{enumerate}
\end{remark}

From Theorem~\ref{theorem2.01} we derive the following theorem. 
\begin{theorem}
\label{theorem2}
Assume that $y(t)\ge 0$ is a continuous on $[0, \infty)$ function,
\begin{equation}
\label{eq20}
\int_0^\infty \omega(y(t))\frac{1}{(1+t)^\alpha}dt <\infty,
\qquad  0< \alpha \le 1,
\end{equation}
\begin{equation}
\label{eq21}
y(t) - y(s) \le \int_s^t f(x,y(x))dx,\qquad 0\le s\le t,
\end{equation}
and there exist constants $a>0$ and $\kappa>0$ such that
\begin{equation}
\label{eq22}
\int_s^t \sup_{0\le \zeta\le a}f(x,\zeta)dx\le \kappa a \frac{t-s}{s^\alpha},
\qquad \kappa>0,\quad t>s\gg 1.
\end{equation}
Then,
\begin{equation}
\label{eq24}
\lim_{t\to\infty}y(t) = 0.
\end{equation}
\end{theorem}

\begin{proof}
Let $\varphi(t):=\frac{1}{(1+t)^\alpha}$, $\alpha \in (0,1]$. 
Then one can easily verify that conditions \eqref{teq20.0} hold 
with $C = 1/2$. 
Condition \eqref{teq22} also holds for this choice of 
$\varphi$ and $\theta = 2\kappa$. 
Thus, Theorem~\ref{theorem2} follows from Theorem~\ref{theorem2.01}.
\end{proof}

\begin{remark} 
{\rm
The assumption $\alpha \in (0,1]$ in \eqref{eq20} is essential: if  
$\alpha >1$, then 
inequality \eqref{teq20.0} does not hold for 
$\varphi(t) = \frac{1}{(1+t)^\alpha}$ 
whatever fixed $C>0$ is.
}
\end{remark}

\begin{corollary}
\label{tcor2}
Let $y(t)\ge 0$ be a continuous function on $[0,\infty)$ and 
\begin{equation}
\label{teq25}
\int_0^\infty \omega\big{(}y(t)\big{)}\varphi(t)dt <\infty,
\end{equation} 
where $\varphi(t)>0$ is a continuous function on $[0,\infty)$. 
Assume that there exists a constant $C>0$ such that
\begin{equation}
\label{teq20.02}
\lim_{t\to\infty} \Big(t - \frac{C}{\varphi(t)}\Big) = \infty,\qquad 
M:=\limsup_{t\to \infty} \frac{\max_{\xi\in[t - \frac{C}{\varphi(t)},t]}\varphi(\xi)}
{\min_{\xi\in[t - \frac{C}{\varphi(t)},t]}\varphi(\xi)}<\infty,
\end{equation}
\begin{equation}
\label{teq26}
y(t) - y(s) \le \int_s^t h(\xi) d\xi,\qquad 0\le s\le t,
\end{equation}
where $h(t)\ge0$, $\forall t\in [0,\infty)$, and
\begin{equation}
\label{teq27}
A:=\limsup_{t\to\infty}\frac{h(t)}{\varphi(t)} <\infty. 
\end{equation} 
Then,
\begin{equation}
\label{teq28}
\lim_{t\to\infty}y(t) = 0.
\end{equation}
\end{corollary}

\begin{proof}
Let
\begin{equation}
\label{teqm1}
f(t,y):= h(t),\qquad t\ge 0,\qquad y\ge 0.
\end{equation}
Let us check that condition \eqref{teq22} holds for this $f(t,y)$ and $a=2A$. 
From \eqref{teqm1} one gets
\begin{equation}
\label{teqm2}
f(t,y) \le 2A\varphi(t),\qquad t\gg 1,\qquad \forall y\ge 0.
\end{equation}
This implies
\begin{equation}
\label{teqm3}
\int_s^t \max_{0\le \zeta\le 2A}f(\xi,\zeta)d\xi \le \int_s^t 2A\varphi(\xi)d\xi
\le 2A(t-s)\varphi(s),\qquad t>s\gg 1. 
\end{equation}
This and Theorem~\ref{theorem2.01} imply \eqref{teq28}. 
\end{proof}

\begin{corollary}
\label{cor2}
Let $y(t)\ge 0$ be a continuous function on $[0,\infty)$, 
\begin{equation}
\label{eq25}
\int_0^\infty \omega\big{(}y(t)\big{)}\frac{1}{(1+t)^\alpha}dt <\infty,
\qquad  0< \alpha \le 1,
\end{equation}
and
\begin{equation}
\label{eq26}
y(t) - y(s) \le \int_s^t h(\xi) d\xi,\qquad 0\le s\le t,
\end{equation}
where $h(t)\ge0$, $\forall t\in [0,\infty)$. 
If 
\begin{equation}
\label{eq27}
A:=\limsup_{t\to\infty}h(t)t^\alpha <\infty, 
\end{equation} 
then
\begin{equation}
\label{eq28}
\lim_{t\to\infty}y(t) = 0.
\end{equation}
\end{corollary}

\begin{proof}
Let $\varphi(t) = \frac{1}{(t+1)^\alpha}$, $\alpha \in (0,1]$.  
Then conditions \eqref{teq20.02} hold with $C=\frac 1 2$ and $M=1$, and 
condition \eqref{teq27} also holds. 
Thus, \eqref{eq28} follows from Corollary~\ref{tcor2}. 
\end{proof}

\begin{theorem}
\label{theorem3}
Assume that  $g\ge 0$ is a continuously differentiable function 
on $[0,\infty)$,
\begin{equation}
\label{eqm29}
\dot{g}(t) \le -a(t)f(g(t)) + b(t),\qquad g(0) = g_0,
\end{equation}
where $f(t)$ is a nonnegative continuous function on $[0,\infty)$, 
$f(0)=0$, $f(t)>0$ for $t>0$, and 
\begin{equation}
\label{eqfff1}
m(\epsilon) := \inf_{x\ge\epsilon} f(x) > 0,\qquad \forall \epsilon>0.
\end{equation}
If $a(t)> 0$, $b(t)\ge 0$ are continuous on $[0,\infty)$ functions, and 
\begin{equation}
\label{eqm30}
\int_0^\infty a(s)ds = \infty,\qquad \lim_{t\to\infty}\frac{b(t)}{a(t)} 
=0,
\end{equation}
then
\begin{equation}
\label{eqm31}
\lim_{t\to\infty} g(t) = 0.
\end{equation}
\end{theorem}

\begin{proof}
Let 
\begin{equation}
\label{eqm32}
s := s(t):=\int_0^t a(\xi)d\xi.
\end{equation}
It follows from \eqref{eqm30} that the map $t\to s$ maps $[0,\infty)$
onto $[0,\infty)$. Let $t(s)$ be the inverse map and 
define $w(s)=g(t(s))$. Then \eqref{eqm29} takes the form
\begin{equation}
\label{eqm33}
w'(s) \le -f(w) + \beta(s),\qquad w(0) = g_0,
\end{equation} 
where
\begin{equation}
\label{eqm34}
w' = \frac{dw}{ds},\qquad \beta(s) = \frac{b(t(s))}{a(t(s))},
\qquad \lim_{s\to\infty} \beta(s) = 0. 
\end{equation}
Assume that \eqref{eqm31} does not hold. 
Then there exist $\epsilon>0$ and $(s_n)_{n=1}^\infty)$ such that
\begin{equation}
\label{eqm35}
0<s_n\nearrow \infty,\qquad w(s_n)>\epsilon,\qquad \forall n.
\end{equation}

From the last relation in \eqref{eqm34} it follows that there exists 
$T>0$ such that
\begin{equation}
\label{eqm36}
\beta(s)< \frac{m(\epsilon)}{2},\qquad \forall s\ge T.
\end{equation}
Since $s_n\nearrow \infty$, there exists $N>0$ such that $s_n>T$, 
$\forall n\ge N$. 
Thus,
\begin{equation}
\label{eqm37}
w'(s_n) \le - f(w(s_n)) + \beta(s_n) \le - m(\epsilon) + 
\frac{m(\epsilon)}{2}<0,\qquad \forall n\ge N. 
\end{equation} 
Since $w(s)$ is continuously differentiable on the interval 
$(s_{n-1},s_{n})$ and $w'(s_n)<0$, $\forall n\ge N$, 
there are two possibilities:

{\it Case 1}: $w'(s) < 0$, $n\ge N$, for all $s\in (s_{n-1},s_{n})$.

{\it Case 2}: there exists a point $t_n\in (s_{n-1},s_{n})$ such that $w'(s)<0$, $\forall s\in(t_n,s_n)$ 
and $w'(t_n)=0$ where $n\ge N$. 

We claim that Case 2 cannot happen if $n\ge N$ is sufficiently large,
namely so large that $\beta(t_n)<m(\epsilon)$. 
Indeed, if Case 2 holds for such $n$, then 
\begin{equation}
\label{eqm38}
w'(t_n) = 0,\qquad w(t_n)>w(s_n)>\epsilon.
\end{equation}
This and \eqref{eqm33} imply
\begin{equation}
\label{eq}
0 = w'(t_n) \le - f(w(t_n)) + \beta(t_n) < -m(\epsilon) + \beta(t_n),
\end{equation}
i.e., $0<m(\epsilon)<\beta(t_n)$.
This contradicts the assumption $\lim_{t\to \infty} 
\beta(t)=0$ because  if $n$ is sufficiently large then 
$t_n$ is so large that $\beta(t_n)< m(\epsilon)$. 

Since Case 2 cannot happen for all sufficiently large $n$,  
there exists $N_1>0$ sufficiently large so that
\begin{equation}
\label{eqfix}
w'(s)<0,\qquad \forall s\in(s_{n-1},s_n),\qquad n\ge N_1.
\end{equation}
Thus,
\begin{equation}
\label{eqfix2}
w'(t)\le 0,\qquad \forall t\ge s_{N_1}.
\end{equation}
Therefore $w(t)$ decays monotonically for all sufficiently large $t$.
Since $w(t)\ge 0$, one concludes that the following limit 
$W\geq 0$ exists and is finite
\begin{equation}
\label{eqfix3}
\lim_{t\to\infty} w(t) = W<\infty.
\end{equation}
This and \eqref{eqm33} imply
\begin{equation}
\label{eqfix4}
\limsup_{t\to\infty}w'(t)\le \lim_{t\to\infty}[-f(w(t)) + \beta(t)] \le 
-m(W).
\end{equation}
If $W\not=0$, then $m(W)>0$ and 
\begin{equation}
\label{eqfix5}
\limsup_{t\to\infty}w'(t) \le -m(W)<0.
\end{equation}
This is impossible since $w(t)\ge 0$, $\forall t$. 
This contradiction 
implies that $W=0$, so  \eqref{eqm31} holds. 

Theorem~\ref{theorem3} is proved. 
\end{proof}

\begin{remark}
Theorem~\ref{theorem3} is proved in \cite{R499}  
under the assumption that $f\in Lip_{loc}[0,\infty)$ and
\begin{equation}
\label{eqrlem}
f(0)=0,\qquad f(u)>0\quad \text{for}\quad t>0,\quad f(u)\ge c>0\quad 
\text{for}\quad u\ge 1,
\end{equation} 
where $c=const$. The assumption $f\in Lip_{loc}[0,\infty)$
was used in \cite{R499} in order to prove the global existence of 
$g(t)$. In this paper we assume the global existence of $g(t)$, and
give a new  simple proof of Theorem~\ref{theorem3}.

\end{remark}

\begin{theorem}
\label{theorem5}
Assume that $g\ge 0$ is a $C^1([0,\infty))-$function, 
\begin{equation}
\label{eqn29}
\dot{g}(t) \le -a(t)f(g(t)) + b(t),\qquad g(0) = g_0,
\end{equation}
where $f(t)\geq 0$ 
is a non-decreasing function on $[0,\infty)$, $f(0)=0$,  $f(t)>0$ if 
$t>0$. If $a(t)>0$, $b(t)\ge 0$ are continuous on $[0,\infty)$ functions, 
and 
\begin{equation}
\label{eqn30}
\int_0^\infty a(s)ds = \infty,\qquad \int_0^\infty \beta(s)ds<\infty,
\qquad \beta(t):=\frac{b(t)}{a(t)},
\end{equation}
then
\begin{equation}
\label{eqn31}
\lim_{t\to\infty} g(t) = 0.
\end{equation}
\end{theorem}

\begin{proof}
Let $s$ be defined in \eqref{eqm32} and $w(s)=g(t(s))$. From \eqref{eqm33} 
one gets 
\begin{equation}
\label{eqn34}
w(t) - w(0) + \int_0^t f(w(s))ds \le \int_0^t \beta(s)ds \le 
\int_0^\infty \beta(s)ds<\infty,\qquad \forall t\ge 0.
\end{equation}
This and the assumption that $w\ge0$ imply 
\begin{equation}
\label{eqn37}
\int_0^\infty f(w(s))ds<\infty. 
\end{equation}
From \eqref{eqm33} one obtains
\begin{equation}
\label{eqn35}
w'(s) \le \beta(s),\qquad \forall s\ge 0.
\end{equation}
Since $\int_0^\infty \beta(s)ds<\infty$, the function $\psi(t):=\int_0^t \beta(s)ds<\infty$ 
is uniformly continuous with respect to $t$ on $[0,\infty)$. 
This, relation \eqref{eqn37}, inequality \eqref{eqn35}, and 
Theorem~\ref{maintheorem} 
imply 
\begin{equation}
\label{eqn38}
\lim_{s\to\infty}w(s) = 0.
\end{equation}
This and the relation $w(s)=g(t(s))$ imply \eqref{eqn31}. 

Theorem~\ref{theorem5} is proved.
\end{proof}

\begin{theorem}
\label{theorem6}
Assume that $g\ge 0$ is a $C^1([0,\infty))-$function,
\begin{equation}
\label{eqn29.2}
\dot{g}(t) \le -a(t)f(g(t)) + b(t),\qquad g(0) = g_0,
\end{equation}
where $f(t)\geq 0$ 
is a non-decreasing continuous function on $[0,\infty)$, $f(0)=0$,  
$f(t)>0$  if $t>0$,  
 $a(t) > 0$ and $b(t)\geq 0$ are continuous functions on $[0,\infty)$, and 
there exists a constant $C>0$ such that  
\begin{equation}
\label{equi5}
\lim_{t\to\infty}\Big(t - \frac{C}{a(t)}\Big) = \infty,\qquad 
\limsup_{t\to\infty} \frac{\max_{\xi\in[t - \frac{C}{a(t)},t]}a(\xi)}{\min_{\xi\in[t - \frac{C}{a(t)},t]}a(\xi)}<\infty.
\end{equation}
If 
\begin{equation}
\label{eqn30.2}
K:=\limsup_{t\to\infty}\frac{b(t)}{a(t)}<\infty,\qquad \int_0^\infty 
b(s)ds<\infty,
%\qquad 
\end{equation}
then
\begin{equation}
\label{eqn31.2}
\lim_{t\to\infty} g(t) = 0.
\end{equation}
\end{theorem}

\begin{proof}
From \eqref{eqn29.2} one gets for all $t\ge 0$ the following inequalities
\begin{equation}
\label{eqg1}
g(t) - g(0) + \int_0^t a(s)f(g(s))ds \le  \int_0^t b(s)ds \le \int_0^\infty b(s)ds <\infty.
\end{equation}
Thus
\begin{equation}
\label{eqg2}
\int_0^\infty a(s)f(g(s))ds <\infty.
\end{equation}
This relation, \eqref{eqn30.2}, \eqref{equi5}, and Corollary~\ref{tcor2} imply \eqref{eqn31.2}.
Theorem~\ref{theorem6} is proved. 
\end{proof}

\begin{remark}
If $a(t) = O(\frac{1}{(1+t)^\gamma})$, $\gamma \in [0,1)$, then conditions \eqref{equi5} hold for 
any $C>0$. If $a(t) = O(\frac{1}{1+t})$ then conditions \eqref{equi5} hold if $C>0$ is sufficiently small. 
\end{remark}

\section{Applications} 

Let $H$ be a real Hilbert space. Consider the following problem
\begin{equation}
\label{eq29}
\dot{u} = A(t,u) + f(t),\qquad u(0) = u_0; \qquad f\in C([0,\infty);H),
\end{equation}
where $u_0\in H$, $A(t,u): [0,\infty)\times H\to H$ is 
continuous with respect to $t$ and $u$.
Assume  that   
\begin{gather}
\label{eq180}
A(t,0) = 0,\qquad \forall t\ge 0,\\
\label{eq41}
\langle A(t,u) - A(t,v), u -v\rangle \le -\gamma(t)\|u-v\|\omega(\|u-v\|),
\qquad u,v\in H,
\end{gather}
where 
$\gamma(t)>0$ for all $t\ge0$ is a continuous function and 
$\omega(t)\ge 0$ is continuous and strictly increasing function 
on $[0,\infty)$, $\omega(0)=0$. 

The above assumptions are standing and are not repeated. 
Assumption \eqref{eq41} means that $A$ is a dissipative operator.
Existence of the solution to problem \eqref{eq29} with such
operators was discussed in the literature (\cite{Mi},\cite{Sh}, 
\cite{R499}).

Let 
$$\beta(t):=\|f(t)\|.$$ 

Consider the following three assumptions:
\begin{itemize}
\item
{Assumption A)
\begin{equation}
\label{eq38}
\int_0^\infty \gamma(t)dt = \infty,\qquad \lim_{t\to\infty}\frac{\beta(t)}{\gamma(t)} = 0.
\end{equation}}

\item{Assumption B)
\begin{equation}
\label{eqiu1}
\int_0^\infty \gamma(t) dt =\infty,\qquad \int_0^\infty \frac{\beta(t)}{\gamma(t)}dt <\infty.  
\end{equation}
}

\item{Assumption C)
\begin{equation}
\label{eq31}
\int_0^\infty \beta(t)dt <\infty,\qquad \gamma(t) = O(\frac{1}{(1+t)^\alpha}), 
\qquad \limsup_{t\to\infty}\beta(t)t^\alpha <\infty ,
\end{equation}
where $\alpha=const \in (0,1]$.
}
\end{itemize}

\begin{remark}{\rm
Assumption \eqref{eq180} is not an essential restriction: if it does
not hold, define $f_1(t):=f(t)+A(t,0)$, and
$A_1(t,u):=A(t,u)-A(t,0)$. Then $A_1(t,u)$ satisfies assumptions
\eqref{eq180} and \eqref{eq41} and $f_1(t)$ plays the role of $f(t)$.
}
\end{remark}

\begin{lemma}
\label{theorem3.1}
If assumptions \eqref{eq180}--\eqref{eq41} hold, then 
there exists a unique global solution $u(t)$ to \eqref{eq29}. 
\end{lemma}

\begin{proof}
{\it Let us first prove the uniqueness of solution to \eqref{eq29}. }

Assume that $u$ and $v$ are two solutions to \eqref{eq29}. Then one gets
\begin{equation}
\label{eq40}
\dot{u} - \dot{v}  = A(t,u) - A(t,v),\qquad t\ge 0.
\end{equation}
Multiply \eqref{eq40} by $u-v$ and use \eqref{eq41} to obtain
\begin{equation}
\label{eq50}
\frac{1}{2}\frac{d}{dt}\|u-v\|^2 = \langle A(t,u) - A(t,v),u-v\rangle 
\le 0.
\end{equation}
Integrating \eqref{eq50} one gets
\begin{equation}
\label{eq51}
\frac{1}{2}\bigg{(}\|u(t)-v(t)\|^2 - \|u(0)-v(0)\|^2\bigg{)}\le 0,\qquad \forall t\ge 0.
\end{equation}
This implies $u(t)=v(t)$, $\forall t\ge 0$, since $u(0)=v(0)$.

{\it Let us prove the local existence of a solution to \eqref{eq29}.}

In this proof an argument similar to the one in \cite{D} or \cite{R499} 
is used. 
Let $u_n(t)$, called Peano's approximation of $u$, solve the following 
equation
\begin{equation}
\label{eq52}
u_n(t) = u_0 + \int_0^t [A(s,u_n(s-\frac{1}{n})) + f(s)]ds, \qquad t\geq 
0,
\end{equation}
and $u_n(t)=u_0,\qquad \forall t\le 0$.

%Let us prove that $u_n(t)$, $n>0$, is uniquely defined for all $t\geq 0$. 
%Assume that the maximal existence interval for $u_n(t)$ is $(-\infty, T]$, 
%where $0\leq T<\infty$. 
%Then the maximal existence interval for $u_n(t-\frac{1}{n})$ 
%is $(-\infty, T+\frac{1}{n}]$, and  equation  \eqref{eq52} 
%implies that the maximal existence interval for
%$u_n(t)$ is $(-\infty, T+\frac{1}{n}]$. This contradicts the definition 
%of $T$. Thus, $T=\infty$, i.e., $u_n(t)$ is defined for all $t\geq 0$.

Fix some positive numbers $r>0$ and $t_1>0$. Let 
$$
B(u_0,r):= \{u:\|u-u_0\|\le r\}, 
$$
and 
\begin{equation}
\label{eq54}
c := \sup_{(t,u)\in [0,t_1]\times B(u_0,r)} (\|A(t,u)\| + 
\|f(t)\|)<\infty. 
\end{equation}
From \eqref{eq52} one gets
\begin{equation}
\label{eq53}
\|u_n(t) - u_0\| \le ct\le r,\qquad  0\le t\le T_1,\qquad T_1:= \min(t_1,\frac{r}{c})>0,
\end{equation}

Define 
\begin{equation}
\label{eq55}
w_{nm}:=u_n(t) - u_m(t),\qquad g_{mn}: = \|w_{mn}(t)\|,\qquad t\ge 0.
\end{equation}
From \eqref{eq52} one obtains
\begin{equation}
\label{eq56}
\begin{split}
g_{mn}(t)\dot{g}_{mn}(t) =& 
\langle A(t,u_n(t-\frac{1}{n})) - A(t,u_m(t-\frac{1}{m})), u_n(t)- u_m(t)\rangle \\
 =&
\langle A(t,u_n(t-\frac{1}{n})) - A(t,u_m(t-\frac{1}{m})), u_n(t-\frac{1}{n})- u_m(t-\frac{1}{m})\rangle \\
 &+ 
\langle A(t,u_n(t-\frac{1}{n})) - A(t,u_m(t-\frac{1}{m})), u_n(t) - u_n(t-\frac{1}{n})\rangle \\
 &+ 
\langle A(t,u_n(t-\frac{1}{n})) - A(t,u_m(t-\frac{1}{m})), u_m(t-\frac{1}{m}) - u_m(t)\rangle .
\end{split}
\end{equation}
From \eqref{eq56}, \eqref{eq54}, and \eqref{eq41}, one obtains
\begin{equation}
\label{eq57}
\begin{split}
\frac{1}{2}\frac{d}{dt} g^2_{nm}(t) 
\le 4c^2(\frac{1}{n}+ \frac{1}{m}),\qquad m,n\ge 0, \quad t\in [0,T_1].
\end{split}
\end{equation}
Integrating \eqref{eq57}, using the relation $g_{mn}(0)=0$, and taking the 
limit as 
$m,n\to \infty$ one obtains
\begin{equation}
\label{eq58}
0\le \lim_{n,m\to\infty} g^2_{nm}(t) \le \lim_{n,m\to\infty} t4c^2(\frac{1}{n}+ \frac{1}{m}) =0,\qquad
\forall t\in [0,T_1].
\end{equation}
It follows from \eqref{eq58} and the Cauchy criterion for 
convergence of a sequence that 
the following limit exists
\begin{equation}
\label{eq59}
u(t): =\lim_{n\to\infty} u_n(t),\qquad 0\le t\le T_1.
\end{equation}
Passing to the limit $n\to \infty$ in equation \eqref{eq52} and using 
the continuity of $A(t,u)$ on $[0,\infty)\times H$ 
and \eqref{eq59} one concludes that $u(t)$ solves the equation
\begin{equation}
\label{eq60}
u(t) = u_0 + \int_0^t [A(s,u(s))+f(s)]ds,\qquad \forall t \in [0,T_1].
\end{equation}
Thus, the local existence of the solution $u(t)$ to equation \eqref{eq29} 
is proved. 

{\it Let us prove the global existence of $u(t)$.}

Assume that $u(t)$ does not exist globally. Let $[0,T]$ be the 
maximal 
existence interval of $u(t)$. Then, $0<T<\infty$. 
It follows from relation \eqref{eq32} that 
\begin{equation}
\label{eqm45}
\|u(t)\| = g(t) < c=const <\infty,\qquad \forall t\in [0,T).
\end{equation}

Let us prove the existence of the finite limit
\begin{equation}
\label{eqm51}
\lim_{t\to T_-} u(t) = u_T.
\end{equation}

Let $z_h(t):= u(t+h) - u(t)$, $0<t\le t+h<T$. From \eqref{eq} one gets 
\begin{equation}
\label{eqm52}
\dot{z}_h(t) = A(t,u(t+h)) - A(t,u(t)) + f(t+h) - f(t),
\end{equation}
where $0<t\le t+h<T$. Multiply \eqref{eqm52} by $z_h(t)$ and get
\begin{equation}
\label{eqm53}
\frac{1}{2} \frac{d}{d t} \|z_h(t)\|^2 \le 
-\gamma(\|z_h(t)\|)\|z_h(t)\|\omega(\|z_h(t)\|) + \|z_h(t)\|\|f(t+h) - f(t)\|.
\end{equation}
This implies
\begin{equation}
\label{eqm54}
\frac{d}{d t} \|z_h(t)\| \le 
 \|f(t+h) - f(t)\|,\qquad 0<t\le t+h<T.
\end{equation}
Integrating \eqref{eqm54} one gets
\begin{equation}
\label{eqm55}
\begin{split}
\|z_h(t)\| &\le \|z_h(0)\| + \int_0^t  \|f(x+h) - f(x)\| dx\\
&\le \|z_h(0)\| + T \max_{0\le t\le T-h} \|f(t+h) - f(t)\|.
\end{split}
\end{equation}
Since $\lim_{h\to +0}\|u(h)-u(0)\|=0$ and 
$\lim_{h\to +0}\max_{0\le t\le T-h} \|f(t+h) - f(t)\|=0$,
one concludes that
\begin{equation}
\label{eqm56}
\lim_{h\to 0}\|u(t+h) - u(t)\|=0,
\end{equation}
and this relation holds uniformly with respect to $t$ and $t+h$ 
such that $t<t+h<T$. Relation \eqref{eqm56} and the Cauchy criterion 
for convergence 
imply the existence of the finite limit in \eqref{eqm51}. 

Consider the following Cauchy problem
\begin{equation}
\label{eqm57}
\dot{u} = A(t,u) + f(t),\qquad u(T) = u_T.
\end{equation}
By the arguments similar to the given above one derives that there 
exists a 
unique solution $u(t)$ to \eqref{eqm57} on $[T,T+\delta]$, where 
$\delta>0$ is a
sufficiently small number.  
From the continuity of 
$f(t)$ and $u(t)$ in $t$ and $A(t,u)$ in both $t$ and $u$ 
one gets
\begin{equation}
\label{eqm58}
\lim_{t\to T_-} \dot{u}(t) = \lim_{t\to T_-} A(t,u) + f(t) 
 = \lim_{t\to T_+} A(t,u) + f(t) = \lim_{t\to T_+} \dot{u}(t),
\end{equation}
and the above limits are finite. 
Thus, the solution to \eqref{eq29} can be extended to the interval 
$[0,T+\delta]$. 
This contradicts the definition of $T$. Thus, $T=\infty$ 
i.e., $u(t)$ exists globally.

Lemma~\ref{theorem3.1} is proved.
\end{proof}

The main result of this Section is the following theorem. 

\begin{theorem}
\label{theoapp}
%Let assumptions \eqref{eq180}--\eqref{eq41} hold and 
If $u(t)$ solves problem \eqref{eq29} globally 
and at least one of the assumptions A), B), or C) holds,  
then  
\begin{equation}
\label{eq32}
\lim_{t\to\infty} \|u(t)\| = 0.
\end{equation}
\end{theorem}

\begin{proof}
The global existence of $u(t)$ follows from Lemma~\ref{theorem3.1}
under the assumptions of this lemma, or from the results in \cite{Mi}
and \cite{Sh}.
 
{\it Let us prove relation \eqref{eq32}.}

Multiplying \eqref{eq29} by $u$, one obtains
\begin{equation}
\label{eq33}
\frac{1}{2} \frac{d}{dt} \|u(t)\|^2 = 
\langle A(t,u),u\rangle + \langle f(t),u\rangle
 \le -\gamma(t)\|u\|\omega(\|u\|) + \|f(t)\| \|u\|.  
\end{equation}
Since $u(t)$ is continuously differentiable, so is 
$ \|u(t)\|$ at the points $t$ at which $\|u(t)\|>0$. 
At these points inequality \eqref{eq33} implies
\begin{equation}
\label{eq34}
\frac{d}{dt} \|u(t)\| \le - \gamma(t)\omega(\|u(t)\|) + \|f(t)\|,\qquad t\ge 0.
\end{equation}
If $ \|u(t)\|=0$ on an open interval  $(a,b)\subset [0,\infty)$,
then $\frac{d}{dt} \|u(t)\|=0$ on $(a,b)$, and inequality
\eqref{eq34} holds trivially, because $\omega(0)=0$ and $\|f(t)\|\geq 0$.
If $ \|u(s)\|=0$ at an isolated point $s>0$, then the right-sided
derivative  of $ \|u(t)\|$ at the point $s$ exists, and
$$\frac{d}{dt} \|u(t)\|=\lim_{\tau \to 
+0}\frac{\|u(s+\tau)\|}{\tau}=\|\frac{d}{dt}u(s)\|,$$
and inequality \eqref{eq34} holds for this derivative.
In what follows we understand by $\frac{d}{dt} \|u(t)\|$ the right-sided
derivative at the points $s$ at which $ \|u(s)\|=0$.
The left-sided derivative of $ \|u(t)\|$ also exists at such points,
and is equal to $-\|\frac{d}{dt}u(s)\|$, but we will not need 
this left-sided derivative.
  
Let $g(t) := \|u(t)\|$ and $\beta(t):=\|f(t)\|$. From \eqref{eq34} one 
gets
\begin{equation}
\label{eq80}
\dot{g}(t) \le -\gamma(t)\omega(g(t)) + \beta(t),\qquad t\ge 0.
\end{equation}
Let $a(t):=\gamma(t)$ and $b(t):=\beta(t)$.

{\it If Assumption A) holds,}
then \eqref{eq32} follows from 
this assumption and Theorem~\ref{theorem3}. 

{\it If Assumption B) holds,} then \eqref{eq32} follows from this assumption and 
Theorem~\ref{theorem5}.  

{\it If Assumption C) holds,}
then \eqref{eq32} follows from this assumption and 
Theorem~\ref{theorem6}. 

Thus, \eqref{eq32} holds. % provided that $u(t)$ exists globally. 

Theorem~\ref{theoapp} is proved.
\end{proof}

{\bf Example.}

Let $D\subset \mathbb{R}^3$ be a bounded domain with a smooth
boundary. Consider problem \eqref{eq29} with
$\gamma(t)=(1+t)^{-\alpha}$, $\alpha=const\in (0,1]$,
$f(t)=O(\frac{1}{(1+t)^k})$, $k>1$, and
$A(t,u):=\gamma(t)[Lu-u^3]$, where  $L$ is a second
order negative-definite selfadjoint Dirichlet elliptic
operator in $D$, e.g.,  $L=\Delta$, where $\Delta$  is the
Dirichlet Laplacian in $D$. Then 
$$\beta(t)=\|f(t)\|\leq \frac c{(1+t)^k},\qquad k>1,$$
conditions \eqref{eq180} and \eqref{eq41}
are satisfied for $u,v\in D(A)$, where $D(A)$ is the domain of definition 
of the operator $A$, $A(u):=\Delta u- u^3$, 
 $D(A)=H^2(D)\cap H^1_0(D)\subset C(D)$, $H^{\ell}$,
$\ell=1,2$, are the usual Sobolev spaces, and 
the inclusion holds by the Sobolev embedding theorem in $R^3$.
The function $\omega (r)$ in \eqref{eq41} in this example is $\omega 
(r)=cr$, where $c>0$ is a constant. This follows from
the known inequality  
$$-\langle Lu,u\rangle=\|\nabla u\|\geq c(D)\|u\|,$$ valid for 
$u\in H^1_0(D)$, $c(D)=const$ does not depend on $u\in H^1_0(D)$. In this 
example the operator $A$ is not continuous in $H$,
but the global solution to  problem \eqref{eq29} exists and is 
unique (see, e.g., \cite{Sh}, \cite{Mi}, \cite{zheng}).
One checks that Assumption C)
is satisfied, and concludes using Theorem 2.14 that 
\eqref{eq32} holds
for the solution to \eqref{eq29} in this example.

Theorem 2.14 can be applied regardless of the method by which 
the global existence of the unique solution to problem (3.1)
is established and inequality (3.31) is derived for this solution.

Let  $\langle \cdot, \cdot\rangle$ denote the inner product and
$\|\cdot\|$ denote the norm in $L^2(D)$.
Then the usual ellipticity constant $c_1=\gamma(t)c(D)$ in the
inequality  
$$c_1\|u\|^2\leq -\gamma (t)\langle Lu,u\rangle$$ 
tends to zero
as $t\to \infty$, so one deals with a degenerate elliptic operator as 
$t\to \infty$ in problem \eqref{eq29} in this example.

One can extend the result in this example to much more general
nonlinearities. For instance, if $A(t,u)=\gamma(t)[Lu-h(u)]$, where 
$uh(u)\geq 0$ for all $u\in R$,  and $h$ satisfies a local Lipschitz
condition, then one can derive an a priori
bound for the solution $u(t)$ of (3.1) $\sup_{t\geq 0}\|u(t)\|\leq c$, 
and prove the global existence and uniqueness of the solution $u(t)$
to problem (3.1) using, for instance, the method from
 \cite{R495}. The assumption $uh(u)\geq 0$ for all $u\in R$
makes it possible to consider nonlinearities $h(u)$ with an arbitrary 
large speed of growth at infinity. Let us outline the derivation of the 
above bound.
Multiplying (3.1) by $u$ and using the estimate 
$\langle Lu,u\rangle\leq -c\|u\|^2$, the assumption $uh(u)\geq 0$,
and denoting $ g:=\|u\|^2$, one gets the following inequality
$$\dot{g}\leq -2c\gamma(t) g + 2\|f\|g^{1/2},\qquad g(0)= \|u_0\|^2.$$
For simplicity and without loss of generality
assume that $u_0=0$. Then  it is not difficult to derive the 
following inequality
$$\|u(t)\|\leq \int_0^t  \|f(s)\|e^{-c\int_s^t\gamma(\tau)d\tau}ds.$$
Using the assumption $k>1$, one obtains from this inequality
the following estimate:
$$\sup_{t\geq 0}\|u(t)\|\leq \frac 1 {k-1}.$$


\begin{thebibliography}{99}

\bibitem{Be} E. F. Beckenbach, R Bellman, {\it Inequalities}, 
Springer-Verlag, 
Berlin, 1961. 

\bibitem{Bu} P.S. Bullen, {\it A dictionary of inequalities}, Longman, 
Essex, 1998. 

\bibitem{D} K. Deimling, {\it Nonlinear functional analysis}, 
Springer-verlag, berlin, 1985.


\bibitem{Mi} I.Miyadera, {\it Nonlinear semigroups}, Amer. Math. Soc., 
Providence RI, 1992.

%\bibitem{KS} P. Krejci and J. Sprekels, 
%Weak stabilization of solutions to PDEs with hysteresis in 
%thermoviscoelastoplasticity, 
%in R.P. Agarwal, F Newman, J Vosmansky eds., EQUASDIFF 9 Proceedings 
%Masaryk Univ., Brno, 1998, 81--96. 

\bibitem{R499} A. G. Ramm, {\it Dynamical Systems Method
for solving operator equations}, 
Elsevier, Amsterdam, 2007.

\bibitem{R495} A. G. Ramm, Existence of the solution to a nonlinear 
equation, J. Math. Anal. Appl., 316, (2006), 764-767. 

\bibitem{Sh} R.E. Showalter, {\it Monotone operators in Banach space and
nonlinear partial differential equations}, Amer. Math.Soc., Providence RI, 
1997.

\bibitem{S} J. E. Slotine and W. Li, {\it Applied nonlinear control}, 
Prentice Hall, New Jersey, 1991.


\bibitem{zheng} S. Zheng, {\it Nonlinear evolution equations}, Chapman \& 
Hall/CRC, 2004. 





%\bibitem{zheng165} S. Zheng, Asymptotic behavior for strong solutions of 
%the Navier-Stokes equations 
%with external forces, {\it Nonlinear Analysis}, 45 (4), (2001), 435-446. 


\end{thebibliography}
\end{document}